\documentclass[a4paper, 11pt, final]{article}

\usepackage[all]{xy}
\usepackage{amssymb, dsfont, amsmath, amsthm}

\usepackage{fullpage}

\newcommand{\num}{\mathrm{num}}
\renewcommand{\hom}{\mathrm{hom}}
\newcommand{\alg}{\mathrm{alg}}
\newcommand{\C}{\mathbf{C}}
\newcommand{\Q}{\mathbf{Q}}
\newcommand{\Z}{\mathbf{Z}}

\newcommand{\im}{\mathrm{Im}\,}

\newcommand{\End}{\mathrm{End}}
\newcommand{\Hom}{\mathrm{Hom}}
\newcommand{\M}{\mathcal{M}}
\newcommand{\id}{\mathrm{id}}

\newcommand{\Spec}{\mathrm{Spec}}

\newcommand{\et}{\text{\'et}}

\renewcommand\M{\mathcal{M}}
\renewcommand\r{\rightarrow}

\newcommand\h{\mathrm{H}}

\renewcommand\h{\mathfrak{h}}

\renewcommand\L{\mathbb{L}}

\newcommand\eff{\mathrm{eff}}
\newcommand\alb{\mathrm{alb}}
\newcommand\Alb{\mathrm{Alb}}

\usepackage[english]{babel}

\newtheorem{theorem}{Theorem}[section]
\newtheorem{lemma}[theorem]{Lemma}
\newtheorem{e-proposition}[theorem]{Proposition}
\newtheorem{corollary}[theorem]{Corollary}
\newtheorem{e-definition}[theorem]{Definition\rm}

\newtheorem{theoreme}{Th\'eor\`eme}[section]

\newtheorem{proposition}[theoreme]{Proposition}

\title{Pure motives with representable Chow groups}

\author{Charles Vial}
\date{}
\begin{document}

\maketitle










\begin{abstract}
Let $k$ be an algebraically closed field. We show using Kahn's and
Sujatha's theory of birational motives that a Chow motive over $k$
whose Chow groups are all representable (in the sense of definition
\ref{def}) belongs to the full and thick subcategory of motives
generated by the twisted motives of curves.

\vskip 0.5\baselineskip

\begin{center}{ \bf R\'esum\'e} \end{center}
\noindent {\bf Motifs purs dont les groupes de Chow sont
  repr\'esentables. } Soit $k$ un corps alg\'ebriquement clos. Nous
prouvons, en nous servant de la th\'eorie des motifs birationnels
d\'evelopp\'ee par Kahn et Sujatha, qu'un motif de Chow d\'efini sur
$k$ dont les groupes de Chow sont tous repr\'esentables (au sens de la
d\'efinition \ref{def}) appartient \`a la sous-cat\'egorie pleine et
\'epaisse des motifs engendr\'ee par les motifs de courbes tordus.
\end{abstract}


\section*{Version fran\c{c}aise abr\'eg\'ee}
Dans cette note, nous pr\'esentons une preuve de l'\'enonc\'e suivant.

\begin{theoreme} [cf. theorem \ref{2}]
  Soit $k$ un corps alg\'ebriquement clos et soit $\Omega \supset k$
  un domaine universel, i.e. un corps alg\'ebriquement clos de degr\'e
  de transcendance infini sur $k$. Soit $M$ un motif de Chow rationnel
  sur $k$ dont les groupes de Chow $CH_j(M_\Omega)_{\alg}$ sont
  repr\'esentables pour tout entier $j$. Alors $M$ est isomorphe \`a
  une somme directe de motifs de Lefschetz et de $\h_1$ de vari\et\'es
  ab\'eliennes tordus.
\end{theoreme}

Ici, $CH_j(M_\Omega)_{\alg}$ d\'esigne le groupe des cycles
alg\'ebriques de dimension $j$ alg\'ebriquement triviaux modulo l'\'equivalence
rationnelle.  La notion de repr\'esentabilit\'e est d\'efinie en
\ref{def}.  Un tel r\'esultat est connu pour les motifs de surfaces et
a \'et\'e r\'ecemment prouv\'e pour les motifs de vari\'et\'es lisses
projectives de dimension $3$ par Gorchinskiy et Guletskii \cite{GG}.
Kimura \cite{Kim} a prouv\'e qu'\'etant donn\'e un motif $M$, si ses
groupes de Chow rationnels $CH_j(M_\Omega)$ sont des $\Q$-espaces
vectoriels de dimension finie alors $M$ est isomorphe \`a une somme
directe de motifs de Lefschetz, offrant ainsi une g\'en\'eralisation
d'un th\'eor\`eme d\^u \`a Jannsen \cite[Th.  3.5.]{Jannsen}. Notre
r\'esultat g\'en\'eralise les r\'esultats cit\'es ci-dessus et redonne
un th\'eor\`eme d\^u \`a Esnault et Levine \cite{EL} qui montre que
pour une vari\'et\'e complexe lisse et projective, si l'application
classe de Deligne rationnelle en tous degr\'es est injective alors
elle est bijective.

Notre m\'ethode repose sur l'existence de projecteurs de Chow relevant
le plus grand facteur direct du motif num\'erique $\bar{M}$ isomorphe
\`a un objet dans la sous-cat\'egorie pleine et \'epaisse des motifs
num\'eriques engendr\'ee par les motifs de courbes tordus. Une telle
construction est l'objet de \cite[\S 1]{Vial3}. Nous nous servons de
fa\c{c}on essentielle de la th\'eorie des motifs purs birationnels
d\'evelopp\'ee par Kahn et Sujatha. Nous aimons penser \`a cette
th\'eorie comme \`a une mani\`ere synth\'etique de proc\'eder \`a une
d\'ecomposition g\'en\'eralis\'ee de la diagonale telle qu'elle a
\'et\'e mise en \oe uvre par Jannsen, et Esnault et Levine entre
autres.

\section{Introduction}
\label{}




Given a field $k$, there are three categories we will be
dealing with.
\begin{itemize}
\item $\M^\eff(k,\Q)$, the category of effective Chow motives with
  coefficients in $\Q$.
\item $\M(k,\Q)$, the category of Chow motives with coefficients in
  $\Q$, see \cite{Scholl}.
\item $\M^\circ(k,\Q)$, the category of birational   Chow motives
  with coefficients in $\Q$, as defined by Kahn and Sujatha in
  \cite{KS2}.
\end{itemize}

Roughly, $\M(k,\Q)$ is obtained from $\M^\eff(k,\Q)$ by inverting the
Lefschetz motive and $\M^\circ(k,\Q)$ is obtained from $\M^\eff(k,\Q)$
by killing the Lefschetz motive (modulo taking the pseudo abelian
envelope).

Objects in $\M(k,\Q)$ are triples $(X,p,n)$ and morphisms are given by
$$\Hom_k((X,p,n),(Y,q,m)) := q \circ CH_{\dim X+n-m}(X \times Y) \circ
p$$ where $CH_i$ denotes the Chow group of $i$-dimensional cycles
tensored with $\Q$. To any smooth projective variety $X$ we associate
functorially the motive $\h(X):=(X,\id_X,0)$. The category
$\M^\eff(k,\Q)$ is the full subcategory of $\M(k,\Q)$ whose objects
have the form $(X,p,0)$.

Let $\mathcal{L}$ be the ideal of $\M^\eff(k,\Q)$ generated by those
morphisms that factor through an object of the form $N \otimes \L$
with $N$ an effective Chow motive and $\L = (\Spec \ k, \id, 1)$ the
Lefschetz motive. Kahn and Sujatha \cite{KS2} define the $\Q$-linear
tensor category $\M^\circ(k,\Q)$ of pure birational Chow motives over
$k$ to be the pseudo-abelianization of the quotient category
$\M^\eff(k,\Q)/\mathcal{L}$. The functor $\M^\eff(k,\Q) \r
\M^\circ(k,\Q)$ will be denoted by $M \mapsto M^\circ$ and to any
smooth projective variety $X$ we associate functorially the motive
$\h^\circ(X)$.

For each of these three categories, we will write $\Hom_k$ for the
groups of morphisms. It will be clear in which category this takes
place. Note that since the functor $\M^\eff(k,\Q) \r \M(k,\Q)$ is
fully faithful, it doesn't matter in which of these two categories we
consider $\Hom_k(M,N)$ for two effective motives $M$ and $N$.

Given a field extension $L/k$, there are base change functors for each
of these three categories. Given a motive (either effective, pure or
birational) $M$ over $k$, we will write $M_L$ for its image in the
corresponding category of motives (either effective, pure or
birational) over $L$. Moreover, for two motives $M$ and $N$ over $k$,
we write
$$\Hom_L(M,N) := \Hom_L(M_L,N_L).$$

An essential feature of Kahn and Sujatha's category of birational
motives is the following (cf. \cite[(2.5)]{KS2})\vspace{2.8pt}

\begin{theorem}[Kahn-Sujatha] \label{morphism} Let $X$ and $Y$ be
  smooth projective varieties over $k$. Denote by $k(X)$ the function
  field of $X$. Then
  $$\Hom_k(\h^\circ(X),\h^\circ(Y)) = CH_0(Y_{k(X)})  =
  \Hom_{k(X)}(\mathds{1}^\circ,\h^\circ(Y)).$$
\end{theorem}

We now fix a field $\Omega$ containing $k$ which is a
universal domain, i.e. an algebraically closed field of infinite
transcendence degree over its prime subfield.  Before we start, we
need to compute some $\Hom$ groups in the category of birational
motives. Let $X$ and $Y$ be smooth projective varieties. Then by
theorem \ref{morphism}
$\Hom_k(\h^\circ(Y),\h^\circ(X)) = \varinjlim CH^{\dim X}(U \times
X)$, where the limit runs through all nonempty open subsets $U$ of
$Y$.  Now assume $M=(X,p)$ is an effective Chow motive and denote by
$M^\circ$ its image in $\M^\circ (k,\Q)$. Then, we have \vspace{-10pt}
\begin{eqnarray} \nonumber \Hom_k ( \h^\circ(Y),M^\circ) & = & p \circ
  \Hom_k (\h^\circ(Y),\h^\circ(X)) = p \circ \varinjlim CH^{\dim X}(U
  \times X) \\ \nonumber & = & \varinjlim \ (\id_U \otimes p)_*
  CH^{\dim X}(U \times X) = (\id_{k(Y)} \otimes p)_* CH_0(X_{k(Y)}) \\
  \nonumber & = & p_{k(Y)} \circ CH_0(X_{k(Y)}) = \Hom_{k(Y)} (
  \mathds{1}^\circ,M^\circ).
\end{eqnarray}

\begin{lemma} \label{b} Let $M \in \M^\eff(k,\Q)$ and let $M^\circ$
  denote its image in $\M^\circ(k,\Q)$. Then the following statements are
  equivalent.
  \begin{itemize}
  \item $M^\circ=0$.
  \item $\Hom_\Omega(\mathds{1}^\circ,M^\circ) = 0$.
  \item There exists an effective Chow motive $N \in \M^\eff(k,\Q)$
    such that $M$ is isomorphic to $N(1)$.
  \end{itemize}
 \end{lemma} \vspace{2.8pt}
 \noindent \textit{Proof.} The first and the third statements are
 equivalent by \cite{KS2}. Moreover, the first statement obviously
 implies the second one. It remains to prove that the second statement
 implies the first one.  Suppose $M = (X,p)$. Then $\Hom_\Omega (
 \mathds{1}^\circ,M^\circ) \supset \Hom_{k(X)} (
 \mathds{1}^\circ,M^\circ) = \Hom_k ( \h^\circ(X),M^\circ) \supset
 \Hom ( M^\circ,M^\circ) = \End(M^\circ).$ Consequently
 $\End(M^\circ)=0$ and thus $M^\circ=0$.  \qed

 \begin{lemma}[see also \cite{GG}, lemma 1] \label{a} Let $M \in
   \M(k,\Q)$.  Then the following statements are equivalent.
  \begin{itemize}
  \item $M=0$.
  \item $\Hom_\Omega(\mathds{1}(i),M) = 0$ for all $i \in \Z$.
  \end{itemize}
\end{lemma}\vspace{2.7pt}

\noindent \textit{Proof.}  Assume $M = (X,p,n)$ is non-zero, $\dim X
=d$ and that $M$ is isomorphic to some effective motive $N = (Y,q,0)$.
Then $\End(M) \simeq \Hom(M,N) \subseteq \Hom(\h(X)(n),N) \simeq \Hom
( \mathds{1}(n+d), \h(X)\otimes N) \subseteq CH_{n +d}(X \times Y)$
and hence $\End(M) \neq 0$ (i.e. $M \neq 0$) implies $n \geq -d$.

Thus there is an integer $j$ which is the smallest integer such that
$M(j)$ is effective.  Then, by assumption,
$\Hom_\Omega(\mathds{1},M(j))=0$. Therefore,
$\Hom_\Omega(\mathds{1}^\circ,M(j)^\circ)=0$. By lemma \ref{b}, this
implies there exists an effective motive $N$ such that $M(j) \simeq
N(1)$. Hence, $M(j-1) \simeq N$ is effective, contradicting the choice
of $j$.  \qed

\section{Representability}

Let $M = (X,p,j)$ be a motive in $\M(k,\Q)$ and let $CH_{i}(M)_\sim :=
\Hom_k (\mathds{1},M(-i))_\sim = p_*CH_{i-j}(X)_{\sim}$ be the
subgroup of $CH_{i}(M)$ made of those cycles $\sim 0$ for an adequate
equivalence relation $\sim$. Such a definition is unambiguous in the
following sense : it can be checked that if $p \in \End (\h(X))$ is an
idempotent then $(p_*CH_i(X))_\sim = p_*(CH_i(X)_\sim)$. In what
follows we will be mainly interested in $\sim = \alg$ where $\alg$
denotes algebraic equivalence. \vspace{2.8pt}

\begin{e-definition} \label{def} Let $\Omega$ be a universal domain
  over $k$. We say that the Chow group $CH_i(M)_{\alg}$ of
  algebraically trivial cycles of a motive $M = (X,p,j) \in \M(k,\Q)$
  is \emph{representable} if there is a smooth projective curve $C$
  over $\Omega$ (not necessarily connected) and a correspondence
  $\Gamma \in \Hom_\Omega (\h_1(C),M(-i))$ such that $ \Gamma_* :
  \Hom_\Omega (\mathds{1},\h_1(C)) \r \Hom_\Omega
  (\mathds{1},M(-i))_\alg$ is surjective.  We say that the total Chow
  group $CH_*(M)_\alg$ of a motive $M \in \M(k,\Q)$ is
  \emph{representable} if $CH_i(M)_\alg$ is representable for all $i$.
\end{e-definition} \vspace{2.8pt}

Notice that we do not require the curve $C$ to be defined over $k$
(which would have been more restrictive). The notion of
representability chosen here seems to be the most appropriate in the
language of motives.  Proposition \ref{repr} below, which generalizes
Jannsen's \cite[1.6]{Jannsen} where idempotents are not being dealt
with, shows that most notions of representability for zero-cycles are
the same.  First we need a lemma whose proof can be found in
\cite[(1.4)--(1.7)]{Saito}.

\begin{lemma} \label{functoriality} Let $X$ and $Y$ be smooth
  projective varieties over an algebraically closed field $F$. Then
  there exists an albanese map $\mathrm{alb}_X : CH_0(X)_\alg \r
  \mathrm{Alb}_X(F)$.  Moreover, if $\alpha \in \Hom_F(\h(X),\h(Y))$
  then $\alpha$ in-duces a homomorphism $\bar{\alpha} :
  \mathrm{Alb}_X(F) \r \mathrm{Alb}_Y(F)$ satisfying $\bar{\alpha}
  \circ \mathrm{alb}_X = \mathrm{alb}_Y \circ \alpha_* : CH_0(X)_\alg
  \r \mathrm{Alb}_Y(F)$.
\end{lemma}

\vspace{2.8pt}

\begin{proposition} \label{repr} Let $M=(X,p) \in \M^\eff(k,\Q)$.
The following statements are equivalent.\vspace{2.8pt}

$i)$ $CH_0(M)_\alg$ is representable.

$i')$ There is a smooth projective curve $C$ over $k$ and a
correspondence $\Gamma \in \Hom_k (\h_1(C),M)$ such that $
(\Gamma_\Omega)_* : CH_0(C_\Omega)_\alg \r
(p_\Omega)_*CH_0(X_\Omega)_\alg$ is surjective.

$i'')$ $(p_\Omega \circ \iota_\Omega)_*CH_0(B_\Omega)_\alg =
(p_\Omega)_*CH_0(X_\Omega)_\alg$ where $\iota : B \hookrightarrow X$
is a smooth linear section of $X$ of dimension $1$.\vspace{2.8pt}

$ii)$ There exists a closed subvariety $Y \subset X_\Omega$ of
dimension $1$ such that for all $\gamma \in CH_0(X_\Omega)$,
$(p_\Omega)_*\gamma$ has vanishing restriction in $CH_0(X_\Omega-Y)$.

$ii')$ There exists a closed subvariety $Y \subset X$ of dimension $1$
such that for all $\gamma \in CH_0(X_\Omega)$, $(p_\Omega)_*\gamma$
has vanishing restriction in $CH_0((X-Y)_\Omega)$.\vspace{2.8pt}

$iii)$ There is a decomposition $p_\Omega = p_1 + p_2$ with $p_1, p_2
\in CH_d(X_\Omega \times X_\Omega)$ such that $p_1$ is supported on
$X_\Omega \times Y$ and $p_2$ is supported on $D \times X_\Omega$ for
some curve $Y \subset X_\Omega$ and some divisor $D \subset X_\Omega$.

$iii')$ There is a decomposition $p = p_1 + p_2$ with $p_1, p_2 \in
CH_d(X \times X)$ such that $p_1$ is supported on $X \times Y$ and
$p_2$ is supported on $D \times X$ for some curve $Y \subset X$ and
some divisor $D \subset X$. \vspace{2.8pt}

$iv)$ The albanese map $\mathrm{alb} : CH_0(X_\Omega)_\alg \r
\Alb_{X_\Omega}(\Omega)$ is injective when restricted to $(p_\Omega)_*
CH_0(X_\Omega)_\alg$.
\end{proposition}
 \vspace{3.5pt}

\noindent \textit{Proof.} Clearly $i'' \Rightarrow i' \Rightarrow i$,
$ii' \Rightarrow ii$ and $iii' \Rightarrow iii$. \vspace{2.8pt}

$i' \Rightarrow ii'$. Let $C$ and $\Gamma$ be as in $(i')$. Clearly,
$(\Gamma_\Omega)_*CH_0(C_\Omega)$ is supported on $Y_\Omega$ for $Y$
the projection on $X$ of a representative of $\Gamma$ on $C \times X$.
Therefore, by localization, this group vanishes in
$CH_0((X-Y)_\Omega)$. \vspace{2.8pt}

$ii' \Rightarrow iii'$. We use Bloch's and Srinivas' technique
\cite{BS}. Let $\tilde{p}$ be the image of $p$ under the natural map
$CH_d(X\times X) \r CH_0( k(X) \times_k X)$. Fix an embedding $k(X)
\subset \Omega$ that extends that of $k$.  The natural map $CH_0(Y_k)
\r CH_0(Y_{L})$ is known to be injective for any smooth variety $Y$
over $k$ and any field extension $L/k$.  This combined to the fact
that $\tilde{p} = p_*\eta_X$ for $\eta_X$ the generic point of $X$
seen as a rational point over $k(X)$ implies that, for $Y$ as in
$ii'$, $\tilde{p}$ has vanishing restriction in $CH_0(k(X) \times_k (X
- Y))$. By the localization exact sequence $\tilde{p}$ is supported on
$k(X) \times_k Y$. Let $p_1$ be an element of $CH_d(X \times_k Y)$
mapping to $\tilde{p}$.  Then $p-p_1$ has vanishing restriction in
$CH_0( k(X) \times_k X)$ and thus, again by the localization exact
sequence, is supported on $D \times X$ for some divisor $D$ on $X$.
Set $p_2 := p -p_1$. \vspace{2.8pt}

$iii' \Rightarrow i'$. Thanks to Chow's lemma on $0$-cycles,
$(p_2)_\Omega$ acts trivially on $CH_0(X_\Omega)$, hence
$((p_1)_\Omega)_*CH_0(X_\Omega)=(p_\Omega)_*CH_0(X_\Omega)$.  If $C$
is the normalization of $Y$ and $\alpha$ is the pullback of $p_1$ in
$CH_d(X \times C)$, we can write $p_1 = \beta \circ \alpha$ with
$\beta \in \Hom_k(\h(C),\h(X))$.
Then we have $((p \circ \beta)_\Omega)_* CH_0(C_\Omega)_\alg \supseteq
((p \circ \beta)_\Omega) \circ (\alpha_\Omega)_* CH_0(X_\Omega)_\alg =
(p_\Omega \circ p_\Omega)_*CH_0(X_\Omega)_\alg = (
p_\Omega)_*CH_0(X_\Omega)_\alg$. We also clearly have $((p \circ
\beta)_\Omega)_* CH_0(C_\Omega)_\alg \subseteq (
p_\Omega)_*CH_0(X_\Omega)_\alg$. Therefore $(i')$ follows by taking the
curve $C$ and the correspondence $\Gamma = p \circ \beta$.\vspace{2.8pt}

By working over $\Omega$ instead of $k$, the exact same arguments
prove $i \Rightarrow ii \Rightarrow iii \Rightarrow i$.\vspace{2.8pt}

$iii \Rightarrow iv$. As for the implication $(iii') \Rightarrow (i')$,
$p_2$ acts trivially on $CH_0(X_\Omega)$ and $p_1$ factors as $\alpha
\circ \beta$ with $\alpha \in \Hom_\Omega(\h(X),\h(C))$ and $\beta \in
\Hom_\Omega(\h(C),\h(X))$ for some smooth projective curve $C$ over
$\Omega$. We have to prove that if $x \in CH_0(X_\Omega)_\alg$
satisfies $\alb (p_*x)= 0$ then $p_*x = 0$. Obviously $\alb (p_*x)= 0$
implies $\bar{\alpha} \circ \alb (p_*x)= 0$. By lemma
\ref{functoriality} $\alb_C ( (\alpha \circ p)_* x)=0$ and because
$\alb_C$ is an isomorphism we get $(\alpha \circ p)_* x=0$. Thus
$(\beta \circ \alpha \circ p)_* x=0$, that is $p_* \circ p_*x=0$ i.e.
$p_*x = 0$.\vspace{2.8pt}

$iv \Rightarrow i''$. Fix $x \in CH_0(X_\Omega)_\alg$. We want to show
that there exists $z \in CH_0(B_\Omega)_\alg$ such that $(p_\Omega
\circ \iota_\Omega)_*z = (p_\Omega)_*x$. It is a fact
\cite[4.3]{Scholl} that the induced map $\Alb_{B_\Omega}(\Omega) \r
\Alb_{X_\Omega}(\Omega)$ is surjective. Thus, using also the
bijectivity of $\alb_{B_\Omega}$, there exists $z \in
CH_0(B_\Omega)_\alg$ such that $\alb((p_\Omega)_*x) =
\bar{\iota}_\Omega \circ \alb(z)$. Now, thanks to lemma
\ref{functoriality}, $\alb((p_\Omega \circ \iota_\Omega)_*z) =
\bar{p}_\Omega(\bar{\iota}_\Omega \circ \alb(z)) = \bar{p}_\Omega(\alb
\circ (p_\Omega)_*x) = \alb((p_\Omega)_*x)$. By assumption on the map
$\alb$ we get $(p_\Omega \circ \iota_\Omega)_*z = (p_\Omega)_*x$.
\qed

\section{Main theorem}

We denote by $\M_0(k,\Q)$ (resp. $\M_1(k,\Q)$) the full thick
subcategory of $\M(k,\Q)$ generated by the $\h_0$'s (resp. $\h_1$'s)
of smooth projective varieties over $k$. Equivalently, $\M_0(k,\Q)$ is
generated by the motives of points and $\M_1(k,\Q)$ is generated by
the $\h_1$'s of curves, see \cite{Scholl}.  We write $\M_\num(k,\Q)$
for the category of motives for numerical equivalence with rational
coefficients. Jannsen \cite{Jannsen3} famously proved that this
category is abelian semi-simple. For a motive $P \in \M(k,\Q)$, let
$\bar{P}$ denote its image in $\M_\num(k,\Q)$. \vspace{2.8pt}

\begin{e-proposition} [see \cite{Vial3} for a proof] \label{prop} Let
  $M$ be an object in $\M_0(k,\Q)$ (resp.  in $\M_1(k,\Q)$). Let $N$
  be any motive in $\M(k,\Q)$. Then any morphism $f : M \r N$ induces
  a splitting $N = N_1 \oplus N_2$ with $N_1$ isomorphic to an object
  in $\M_0(k,\Q)$ (resp. in $\M_1(k,\Q)$) and $\bar{N}_1 \simeq \im
  \bar{f}$.
\end{e-proposition} \vspace{2.8pt}

From now on, $k$ is an algebraically closed field and $\Omega$ denotes
a universal domain over $k$. Recall that if $M$ and $N$ are two
Chow motives over $k$, then $\Hom_k(\bar{M},\bar{N}) =
\Hom_\Omega(\bar{M},\bar{N})$.\vspace{2.8pt}

\begin{lemma} \label{cut} Let $M \in \M^\eff(k,\Q)$ and let $n$ be the
  dimension of the finite dimensional vector space
  $\Hom_\Omega(\bar{\mathds{1}},\bar{M})$. Then $\mathds{1}^{\oplus
    n}$ is a direct summand of $M$.
\end{lemma}\vspace{2.8pt}
\noindent \textit{Proof.}  Pick a basis $(\bar{e}_i)_{1\leq i \leq n}$
of the group
$\Hom_k(\bar{\mathds{1}},\bar{M})=\Hom_\Omega(\bar{\mathds{1}},\bar{M})$
of $0$-cycles modulo numerical equivalence on $M$. Lift it to a family
$({e}_i)_{1\leq i \leq n}$ of the Chow group $\Hom_k(\mathds{1},M)$
and consider the morphism $\oplus e_i: \mathds{1}^{\oplus n}
\stackrel{}{\longrightarrow} M.$ By proposition \ref{prop}, $M$ has
then a direct summand $N$ isomorphic to an object in $\M_0$ whose
reduction modulo numerical equivalence is $\bar{\mathds{1}}^{\oplus
  n}$.  Therefore $N \simeq \mathds{1}^{\oplus n}$.  \qed
\vspace{2.8pt}

\begin{lemma} \label{representable} Let $M = (X,p) \in \M^\eff(k,\Q)$
  be such that $ \Hom_k (\mathds{1},M)_\alg$ is representable. Assume
  moreover that $\bar{M}$ has no direct factor of
  the form $\bar{\h}_1(J)$ for an abelian variety $J$. Then
  $\Hom_\Omega (\mathds{1},M)_\alg =0$.
\end{lemma} \vspace{2.8pt}

\noindent \textit{Proof.}  Thanks to proposition \ref{repr} and its
proof (specifically the statement $(i) \Rightarrow (iii')$ plus an extra
argument included in the proof of $(iii') \Rightarrow (i')$) the
representability assumption on $ \Hom_k (\mathds{1},M)_\alg$ yields a
decomposition $p = p_1 + p_2 \in CH_d(X \times X)$ such that $p_1$
factors through a smooth projective curve $C$ over $k$ and $p_2$ is
supported on $D \times X$ for some proper subscheme $D$ of $X$. In
particular $(p_2)_\Omega$ acts trivially on $0$-cycles on $X_\Omega$.
Let's write $p_1 = \beta \circ \alpha$ with $\alpha \in
\Hom_k(\h(X),\h(C))$ and $\beta \in \Hom_k(\h(C),\h(X))$.  The
correspondence $(p_1)_\Omega$ acts as the identity on $ \Hom_\Omega
(\mathds{1},M)$. If $\pi_1$ denotes the projector on $\h_1(C)$ with
respect to the choice of a $0$-cycle of degree $1$ on $C$ (see e.g.
\cite{Scholl}) and if $q_1 := p \circ \beta \circ \pi_1 \circ \alpha$
then $(q_1)_\Omega$ acts as the identity on $ \Hom_\Omega
(\mathds{1},M)_\alg$.  Therefore $(q_1 \circ q_1)_\Omega$ also acts as
the identity on $ \Hom_\Omega (\mathds{1},M)_\alg$. By assumption on
$M$, the map $p \circ \beta \circ \pi_1$ must be numerically trivial.
Hence the map $\pi_1 \circ \alpha \circ p \circ \beta \circ \pi_1 \in
\End(\h_1(C))$ is also numerically trivial. Because $\End(\h_1(C)) =
\End(\bar{\h}_1(C))$ we get that $\pi_1 \circ \alpha \circ p
\circ \beta \circ \pi_1 =0$ and therefore that $q_1 \circ q_1 = p
\circ \beta \circ \pi_1 \circ \alpha \circ p \circ \beta \circ \pi_1
\circ \alpha= 0$. This proves that $ \Hom_\Omega (\mathds{1},M)_\alg =
0$.  \qed \vspace{2.8pt}

\begin{theorem} \label{2} Let $M \in \M(k,\Q)$. Then $CH_*(M)_\alg$ is
  representable if and only if $M$ is isomorphic to a sum of Lefschetz
  motives and twisted $\h_1$'s of abelian varieties.
\end{theorem} \vspace{2.8pt}

\noindent \textit{Proof.}  Assume $M=(X,p,n)$ with $X$ a smooth
projective variety over $k$. Up to tensoring with $\mathds{1}(-n)$ we
can assume that $M$ is effective.  The integers $r$ for which
$CH_r(M_\Omega) := \Hom_\Omega (\mathds{1}(r),M)$ is possibly non-zero
are non-negative.  We proceed by induction on $\mu(M) := \max\{r \ : \
\Hom_\Omega (\mathds{1}(r),M) \neq 0\} \in \{-\infty\}\cup \Z_{\geq
  0}.$

In the case $\mu(M) = -\infty$, that is by definition in the case when
$CH_r(M_\Omega) = 0$ for all integers $r$, we conclude directly by
lemma \ref{a} that $M=0$. Let then $M$ be an effective motive with
$\mu(M) > - \infty$.  Let $n$ be the dimension of the $\Q$-vector
space $\Hom_\Omega(\bar{\mathds{1}},\bar{M})$. By lemma \ref{cut},
there exists a motive $M'$ over $k$ such that $M = \mathds{1}^n \oplus
M'$ and $\Hom_\Omega(\bar{\mathds{1}},\bar{M}') = 0.$

Let $C$ be a curve over $k$ and $\Gamma \in \Hom_k (\h_1(C),M)$ be
such that $\bar{\Gamma} \in \Hom_k (\bar{\h}_1(C),\bar{M}')$ has
maximal image inside $\bar{M}'$ among all curves $C'$ and all
morphisms in $\Hom_k (\bar{\h}_1(C'),\bar{M}')$.
By proposition \ref{prop}, $\Gamma$ induces a splitting $M' = \h_1(J)
\oplus N$ for some abelian variety $J$ and some effective motive $N$
satisfying $\Hom_k (\bar{\h}_1(C'),\bar{N})=0$ for all curves $C'$.
Since $N$ is a direct summand of $M$, the group $\Hom_k
(\mathds{1},N)_\alg$ is representable. By lemma \ref{representable},
$\Hom_\Omega (\mathds{1},N)_\alg = 0$. Moreover, because $N$ is a
direct summand of $M'$, we have $\Hom_\Omega
(\bar{\mathds{1}},\bar{N}) = 0$.  Algebraic equivalence and numerical
equivalence agree on $0$-cycles.  Therefore we have a decomposition
$M= \mathds{1}^{\oplus n} \oplus \h_1(J) \oplus N$ with $\Hom_\Omega
(\mathds{1},N) = 0$.  Hence $\Hom_\Omega (\mathds{1}^\circ,N^\circ) =
0$.  Therefore, by lemma \ref{b}, there exists an effective Chow
motive $N' \in \M^\eff(k,\Q)$ such that $M = \mathds{1}^{\oplus n}
\oplus \h_1(J) \oplus N'(1).$ The Chow group of the motive $N'_\Omega$
is a subgroup of the Chow group of $M_\Omega$, it is therefore
representable.  Clearly $\mu(N') \leq \mu(M) -1$ which concludes the
proof by induction.  \qed \vspace{3.5pt}

\begin{corollary}[Kimura \cite{Kim}] \label{1} Let $M \in \M(k,\Q)$.
  Then $CH_*(M_\Omega)$ is a finite dimensional $\Q$-vector space if
  and only if $M$ is isomorphic to a sum of Lefschetz motives.
\end{corollary} \vspace{2.8pt}

Let $X$ be a smooth projective variety of dimension $d$ over an
algebraically closed subfield $k$ of $\C$. The Abel-Jacobi map $AJ_i :
CH_i(X_\C)_\hom \r J_i(X_\C)$ defined by Griffiths (here Chow groups
are not tensored with $\Q$) restricts to $ CH_i(X_\C)_\alg$ and the
image of the composite map $CH_i(X)_\alg \r CH_i(X_\C)_\alg \r
J_i(X_\C)$ defines an abelian variety over $k$ that we denote
$J_i^a(X)$. \vspace{2.8pt}

\begin{corollary} \label{rep} Assume that the total Chow group of $X$
  is representable. Then,\vspace{-1pt} $$\h(X)= \mathds{1} \oplus
  \h_1(\mathrm{Alb}_X) \oplus \L^{\oplus b_{2}} \oplus
  \h_1(J^a_1(X))(1) \oplus (\L^{\otimes 2})^{\oplus b_{4}} \oplus
  \ldots \oplus \h_1(J^a_{d-1}(X))(d-1) \oplus \L^{\otimes d}$$ where
  $b_i$ denotes the $i^\mathrm{th}$ Betti number of $X$. Moreover,
  algebraic equivalence agrees with numerical equivalence on $X$ and
  the generalized Hodge conjecture holds for $X$.
\end{corollary} \vspace{2.8pt}

\begin{corollary} [Esnault-Levine \cite{EL}]
  Let $X$ be a complex smooth projective variety. Suppose that the
  total rational Deligne cycle class map $cl_\mathcal{D} : \oplus_i
  CH^i(X) \r \oplus_i H_\mathcal{D}^{2i}(X,\Q(i))$ is injective. Then
  it is surjective.
\end{corollary}

\section*{Acknowledgements} Thanks to Magdalene College, Cambridge and
to the EPSRC for financial support under the grant EP/H028870/1. The
idea of using birational motives grew out of a discussion with Bruno
Kahn at the workshop ``Finiteness results for motives and motivic
cohomology'' in Regensburg in February 2009. I am grateful to Bruno
Kahn for explaining to me lemma \ref{b} and related ideas, as well as
to Uwe Jannsen for organizing the workshop and giving me the
opportunity to present my work. I also wish to thank the referee for
his very useful and numerous comments, and also for providing a
reference to lemma \ref{functoriality}.

\end{document}